\documentclass[12pt]{article}
\usepackage{mathrsfs}
\usepackage{pifont}

\usepackage{amsfonts}
\usepackage{latexsym}
\usepackage{amsmath}
\usepackage{amssymb}
\usepackage{color}

\def\ep{\varepsilon}
\def\a{\alpha}
\def\s{\sigma}
\def\g{\gamma}
\def\p{\varphi}
\def\b{\beta}

\title{\textbf   The viability property of jump diffusion processes on Riemannian manifolds}

\date{}

\author{\  Xuehong Zhu \\
 {School of Science, \ Nanjing
University of Aeronautics and
Astronautics }\\
 Nanjing, \ 210016, \ China\\
{E-mail: hilda2002@163.com }}

\begin{document}

\maketitle

\begin{abstract}
In this note, we consider the necessary and sufficient condition for
viability property of diffusion processes
  with jumps on closed submanifolds of $R^{m}$ with some concrete examples.\\
\par  $\textit{Keywords:  Viability property; Viscosity solution; Riemannian
manifold.}$
\end{abstract}



{\section*{\large \bf1. Introduction}}

The viability property had been widely studied in the deterministic
case and a little bit less in the stochastic case. The major
contributions of 1970's and 1980's are quoted in Aubin and Da Prato
[1] and Gautier and Thibault [3]. But these papers were all based on
the stochastic tangent cone and they are all just sufficient
conditions. R. Buckdahn, S. Peng, M. Quincampoix and C. Rainer [2]
used a new method to get the necessary and sufficient condition for
the viability property of SDEs with control. They related viability
with a kind of optimal control problem and took advantage of
comparison theorem of viscosity solutions to some H-J-B equation. S.
Peng and X. Zhu [8] generalized the result to the case with jumps.
M. Michta [7], L. Mazliak [5] and L. Mazliak and  C. Rainer [6] were
interested in formulating weak notions of viability which may be
satisfied more easily.

To be viable on a close submanifold $K\subset R^{m}$, Elton P. Hsu
[4] considered the processes driven by continuous semimartingales in
the form of stratonovich integral. He showed that when the smooth
vector fields are tangent to $K$ along $K$, the process won't leave
$K$ before its explosion time. The approach based on the distance
function of $K$ and Gronwall's inequality, and the result is just
sufficient.

In this paper, we consider SDEs driven by a Brownian motion and a
Poisson process. As an application of the result in [8], we study
the necessary and sufficient condition for which the solution can be
viable on closed submanifolds of $R^{m}$. \vspace{5mm}

 Let $(\Omega,{\cal{F}},P,({\cal{F}}_{t})_{t\geq 0})$ be a complete
stochastic basis such that $\mathcal{F}_{0}$ contains all $P$-null
elements of ${\cal{F}}$, and
$\mathcal{F}_{t^{+}}:=\cap_{\ep>0}\mathcal{F}_{t+\ep}=\mathcal{F}_{t},t\geq
0$, and $\mathcal{F}=\mathcal{F}_{T}$, and suppose that the
filtration is generated by the following
two mutually independent processes:\\
(i) a $d$-dimensional standard Brownian motion $(W_{t})_{0\leq t\leq
T}$, and\\
(ii) a stationary Poisson random measure $N$ on
$(0,T]\times E$, where $E\subset R^{l}\setminus\{0\}$, $E$ is
equipped with its Borel field $\mathscr{B}_{E}$, with compensator
$\hat{N}(dtde)=dtn(de)$, such that $n(E)<\infty$, and
$\{\tilde{N}((0,t]\times A)=(N-\hat{N})((0,t]\times A)\}_{0<t\leq
T}$ is an $\mathcal{F}_{t}$-martingale, for each $A\in
\mathscr{B}_{E}$.

By $T>0$ we denote the finite real time horizon.

We consider a jump diffusion process as follows:
$$
X^{t,x}_{s}=x+\int^{s}_{t}b(r,X^{t,x}_{r})dr+\int^{s}_{t}\sigma(r,X^{t,x}_{r})dW_{r}\\
+\int^{s}_{t}\int_{E}\gamma(r,X^{t,x}_{r-},e)\tilde{N}(drde),s\in
[t,T], \eqno(1.1)
$$
where
$$
\begin{array}{l}
b:[0,\infty)\times R^{m} \rightarrow R^{m},
\gamma:[0,\infty)\times R^{m} \times R^{l} \rightarrow R^{m},\\
\sigma=\{\sigma_{\alpha}^{i}\}:[0,\infty)\times R^{m} \rightarrow
R^{m\times d},i=1,2,...,m,\a=1,2,...,d.
\end{array}
$$

{\bf Definition 1.1.} The SDE (1.1) enjoys the stochastic viability
property (SVP in short) in a given closed set  $K\subset R^{m}$ if
and only if: for any fixed time interval $[0,T]$, for each $(t,x)\in
[0,T]\times K$, there exists a probability space
 $(\Omega,{\cal{F}},P)$, a $d-$dimensional Brownian motion $W$, a
 stationary Poisson process $N$,  such that
 $$
 X^{t,x}_{s}\in K ,\mbox{ \ } \forall \mbox{ \ }s  \in [t,T] \mbox{ \
 } \mbox{ \ }P-a.s..
 $$

We assume that, there exists a sufficiently large constant $\mu>0$
and a function $\rho: R^{l}\rightarrow R_{+}$ with
$$
\int_{E}\rho^{2}(e)n(de)<\infty,
$$
such that

(A1)$b,\sigma,\gamma\mbox{ \ }\mbox{are continuous in}\mbox{ \ }
(t,x),$

(A2) for all $x,x'\in R^{m}$, $t\in [0,+\infty)$
$$
\begin{array}{ll}
| b(t,x)-b(t,x')|+|\s(t,x)-\s(t,x')|\leq \mu|x-x'|,\\
|b(t,x)|+|\s(t,x)|\leq \mu(1+|x|),\\
|\g(t,x,e)-\g(t,x',e)|\leq \rho(e)|x-x'|,\forall e\in E,\\
|\g(t,x,e)|\leq \rho(e)(1+|x|),\forall e\in E.
\end{array}
$$
Here $\langle\cdot\rangle$ and $|\cdot|$ denote, respectively, the
Euclidian scalar product and norm. Obviously under the above
assumptions there exists a unique strong solution to SDE (1.1). We
set $C$ is a constant such that
$$
C\geq 1+2\mu+\mu^{2}+\int_{E}\rho^{2}(e)n(de).
$$

We denote by $C_{2}([0,T]\times R^{m})$ (resp,.
$C^{1,2}_{2}([0,T]\times R^{m})$) the set of all functions in
$C([0,T]\times R^{m})$ (resp., $C^{1,2}([0,T]\times R^{m})$)with
quadratic growth in $x$. In fact, the SVP in $K$ is related to the
following PDE:
$$
 \left\{
\begin{array}{l}
\mathscr{L}u(t,x)+\mathscr{B}u(t,x)-Cu(t,x)+d^{2}_{K}=0,\mbox{ \
}(t,x)\in (0,T)\times R^{m},\\
u(T,x)=d^{2}_{K}(x),
\end{array}
\right. \eqno{(1.2)}
$$
where we denote, for $\p\in C_{2}^{1,2}([0,T]\times R^{m})$ ,
$$
\mathscr{L}\p(t,x):=\frac{\partial \p(t,x)}{\partial t}+\langle
D\p(t,x),b(t,x)\rangle +\frac{1}{2}tr[D^{2}\p(t,x)\s\s^{T}(t,x)],
$$
$$
\mathscr{B}\p(t,x):=\int_{E}[\p(t,x+\g(t,x,e))-\p(t,x)-\langle
D\p(t,x),\g(t,x,e)\rangle]n(de).
$$

{\bf Definition 1.2.} We say a function $u\in C_{2}([0,T]\times
R^{m})$ is a viscosity supersolution (resp., subsolution) of (1.2)
if, $u(T,x)\geq d^{2}_{K}(x)$ (resp., $u(T,x)\leq d^{2}_{K}(x)$) and
for any $\p\in C_{2}^{1,2}([0,T]\times R^{m})$ and any point
$(t,x)\in [0,T]\times R^{m}$ at which $u-\p$ attains its minimum
(resp., maximum),
$$
\mathscr{L}\p(t,x)+\mathscr{B}\p(t,x)-C\p(t,x)+d^{2}_{K}\leq
0,\mbox{ \ (resp.,  }\geq 0).
$$
u is called a viscosity solution if it is both viscosity
supersolution and subsolution.

Now let us recall the characterization of SVP of SDE (1.1) in $K$
(see [8]):

{\bf Lemma 1.3.} We assume (A1) and (A2). Then the following claims
are
equivalent:\\
{\rm(i)}SDE (1.1) enjoys the SVP in $K$; \\
{\rm(ii)}$d^{2}_{K}(\cdot)$ is a viscosity supersolution of PDE
(1.2).

As the end of this section, we define for any closed set $K\subset
R^{m}$ the projection of a point $a$ onto $K$:
$$
\Pi_{K}(a):=\{b\in K|\mbox{ \ } \|a-b\|=\min\limits_{c\in
K}\|a-c\|=d_{K}(a)\}.
$$

In the next section, we will use Lemma 1.3 to get the necessary and
sufficient condition for the viability property of (1.1) in closed
submanifolds of $R^{m}$. Finally, a special case: $S^{2}$ is
considered with some concrete examples.

\vspace{3mm}

{\section*{\large \bf2. Viability in closed submanifolds of
$R^{m}$}}

If $K$ is a closed submanifold of $R^{m}$ without boundary, the
function $d_{K}^{2}(x) $ is smooth in a neighborhood of $K$ (See
[4]). We assume that $V_{\sigma}$ is a smooth vector field on
$R^{m}$ and $\forall f\in C^{2}(R^{m})$, $V_{\sigma}f(x)=\langle
Df(x),\sigma(x)\rangle$. From [4], we have the following lemma:

{\bf Lemma 2.1.} If the vector field $V_{\sigma}$ is tangent to $K$
along $K$, in a sufficiently small neighborhood $U$ of $K$, there
exists $C^{*}$ depending on $|x|(x \in U)$ such that
$$
|V_{\sigma}d_{K}^{2}(x)|\leq
C^{*}d_{K}^{2}(x),|V_{\sigma}V_{\sigma}d_{K}^{2}(x)|\leq
C^{*}d_{K}^{2}(x).
$$

Then we can have the following main result:

{\bf Theorem 2.2. }Under the assumptions (A1) and (A2), SDE (1.1)
enjoys SVP in a closed submanifold $K$ if and only if: $\forall t\in
[0,T], \bar{x}\in K$,
$$
\left\{\begin{array}{l} 2\langle
b(t,\bar{x}),m(\bar{x})\rangle-\sum\limits_{i=1}\limits^{d}\langle\langle
D\sigma_{\a},\sigma_{\a}\rangle(t,\bar{x}),m(\bar{x})\rangle-2\displaystyle\int_{E}\langle\gamma(t,\bar{x},e),m(\bar{x})\rangle
n(de)=0,\\
\langle\sigma_{\a}(t,\bar{x}),m(\bar{x})\rangle=0,\forall
\a=1,2,...,d,\\
\bar{x}+\gamma(t,\bar{x},e)\in K,n(de)-a.s.,
\end{array} \right.
\eqno(2.1)
$$
where $m(\bar{x})$ is any normal vector of $K$ at $\bar{x}$ and
$$\langle D\sigma_{\a},\sigma_{\a}\rangle:= (\langle
D\sigma^{1}_{\a},\sigma_{\a}\rangle,\langle
D\sigma^{2}_{\a},\sigma_{\a}\rangle,...,\langle
D\sigma^{m}_{\a},\sigma_{\a}\rangle )'.$$

{\bf Proof: }All we have to do is to prove that (2.1) is equivalent
to the following statement: $d_{K}^{2}(x)$ is a viscosity
supersolution of PDE (1.2).

(a) If $d_{K}^{2}(x)$ is a viscosity supersolution of PDE (1.2), in
a sufficiently small neighborhood $U$ of $K$ we have: $\forall t \in
[0,T],$
$$
\begin{array}{ll}
&2\langle
b(t,x),x-\Pi_{K}(x)\rangle+\frac{1}{2}tr[D^{2}d_{K}^{2}(x)\sigma\sigma'(t,x)]\\
&+\displaystyle\int_{E}[d_{K}^{2}(x+\gamma(t,x,e))-d_{K}^{2}(x)-2\langle
\gamma(t,x,e),x-\Pi_{K}(x)\rangle]n(de)\\
\leq &(C-1)|x-\Pi_{K}(x)|^{2}.
\end{array}
\eqno(2.2)
$$
Select $x=\bar{x}\in K$ in (2.2), we have
$$
\frac{1}{2}tr[D^{2}d_{K}^{2}(\bar{x})\sigma\sigma'(t,\bar{x})]
+\int_{E}d_{K}^{2}(\bar{x}+\gamma(t,\bar{x},e))n(de)\leq 0.
$$
Since $D^{2}d_{K}^{2}(\bar{x})\geq 0$, we have
$$\bar{x}+\gamma(t,\bar{x},e)\in K, n(de)-a.s.,\eqno{(2.3)}$$
and for all $\a=1,2,...,d$,
$$
\frac{1}{2}V_{\sigma_{\a}}V_{\sigma_{\a}}d_{K}^{2}(x)=\{\langle\langle
D\sigma_{\a},\sigma_{\a}\rangle(t,x),x-\Pi_{K}(x)\rangle+
\frac{1}{2}tr[D^{2}d_{K}^{2}(x)\sigma_{\a}(\sigma_{\a})'(t,x)]\}|_{x\in
K}=0.
$$
On the other hand,
$$
V_{\sigma_{\a}}V_{\sigma_{\a}}d_{K}^{2}(x)=\{2\langle
Dd_K(x),\sigma_{\a}(t,x)\rangle^{2}+2d_K(x)V_{\sigma_{\a}}\langle
Dd_K,\sigma_{\a}\rangle(t,x)\}.
$$
Therefore $$\langle\sigma_{\a}(t,x),Dd_K(x)\rangle|_{x\in
K}=0,\forall \a=1,2,...d.$$ Notice that when $x\notin K$,
$$Dd_K(x)=\frac{x-\Pi_{K}(x)}{|x-\Pi_{K}(x)|},$$ so
$$\langle\sigma_{\a}(t,\bar{x}),m(\bar{x})\rangle=0,\forall
\a=1,2,...d.$$

Divide (2.2) by $|x-\Pi_{K}(x)|$ when $x \notin K$, we get
$$
\begin{array}{ll}
&2\langle
b(t,x),Dd_K(x)\rangle-\sum\limits_{\a=1}\limits^{d}\langle\langle
D\sigma_{\a},\sigma_{\a}\rangle(t,x),Dd_K(x)\rangle+\sum\limits_{\a=1}\limits^{d}\frac{1}{2|x-\Pi_{K}(x)|}V_{\sigma_{\a}}V_{\sigma_{\a}}d_{K}^{2}(x)\\
&+\displaystyle\int_{E}[\frac{d_{K}^{2}(x+\gamma(t,x,e))}{|x-\Pi_{K}(x)|}-|x-\Pi_{K}(x)|-2\langle
\gamma(t,x,e),Dd_K(x)\rangle]n(de)\\
\leq &(C-1)|x-\Pi_{K}(x)|.
\end{array}
$$
By Lemma 2.1, (2.3) and the Lipschtz condition of $\gamma$ in $x$,
when $x\rightarrow \Pi_{K}(x)$ along $x-\Pi_{K}(x)$, we have
$$
\sum\limits_{\a=1}\limits^{d}\frac{1}{2|x-\Pi_{K}(x)|}V_{\sigma_{\a}}V_{\sigma_{\a}}d_{K}^{2}(x)
+\int_{E}\frac{d_{K}^{2}(x+\gamma(t,x,e))}{|x-\Pi_{K}(x)|}n(de)\rightarrow
0.
$$
So let $x\rightarrow \Pi_{K}(x)$  we get
$$
2\langle b(t,\bar{x}),\pm
m(\bar{x})\rangle-\sum\limits_{\a=1}\limits^{d}\langle\langle
D\sigma_{\a},\sigma_{\a}\rangle(t,\bar{x}),\pm
m(\bar{x})\rangle-2\int_{E}\langle\gamma(t,\bar{x},e),\pm
m(\bar{x})\rangle n(de)\leq 0.
$$
So
$$
2\langle
b(t,\bar{x}),m(\bar{x})\rangle-\sum\limits_{\a=1}\limits^{d}\langle\langle
D\sigma_{\a},\sigma_{\a}\rangle(t,\bar{x}),
m(\bar{x})\rangle-2\int_{E}\langle\gamma(t,\bar{x},e),m(\bar{x})\rangle
n(de)=0.
$$

(b) If (2.1) is true, we want to show $d_{K}^{2}(x)$ is a viscosity
supersolution of (1.2). In fact, we just need to prove that $\forall
t \in (0,T)$, (2.2) is true at some neighborhood $U$ of $K$ (See
[9]). So due to (2.1), Lemma 2.1 and the assumption (A2), we have,
$\forall t\in (0,T),x\in U$,
$$
\begin{array}{ll}
&2\langle
b(t,x),x-\Pi_{K}(x)\rangle+\frac{1}{2}tr[D^{2}d_{K}^{2}(x)\sigma\sigma'(t,x)]\\
&+\displaystyle\int_{E}[d_{K}^{2}(x+\gamma(t,x,e))-d_{K}^{2}(x)-2\langle
\gamma(t,x,e),x-\Pi_{K}(x)\rangle]n(de)\\
\leq & 2\langle b(t,x)-b(t,\Pi_{K}(x)),x-\Pi_{K}(x)\rangle+2\langle b(t,\Pi_{K}(x)),x-\Pi_{K}(x)\rangle\\
&-\sum\limits_{\a=1}\limits^{d}\langle\langle
D\sigma_{\a},\sigma_{\a}\rangle(t,\Pi_{K}(x)),x-\Pi_{K}(x)\rangle-2\displaystyle\int_{E}\langle\gamma(t,\Pi_{K}(x),e),x-\Pi_{K}(x)\rangle
n(de)\\
&+\sum\limits_{\a=1}\limits^{d}\langle\langle
D\sigma_{\a},\sigma_{\a}\rangle(t,\Pi_{K}(x))-\langle
D\sigma_{\a},\sigma_{\a}\rangle(t,x),x-\Pi_{K}(x)\rangle\\
 &+\sum\limits_{\a=1}\limits^{d}\langle\langle
D\sigma_{\a},\sigma_{\a}\rangle(t,x),x-\Pi_{K}(x)\rangle+
\frac{1}{2}tr[D^{2}d_{K}^{2}(x)\sigma\sigma'(t,x)]\\
&+\displaystyle\int_{E}[(x+\gamma(t,x,e)-\Pi_{K}(x)-\gamma(t,\Pi_{K}(x),e))^{2}
-d_K^2(x)\\
& \ \ \ \ \ \ \ -2\langle\gamma(t,x,e)-\gamma(t,\Pi_{K}(x),e),x-\Pi_{K}(x)\rangle]n(de)\\
\leq &(C-1)|x-\Pi_k(x)|^2,
\end{array}
$$
where $C$ is a constant depending on $|x|(x\in U)$ and can be chosen
large enough such that
$$C \geq 1+2\mu+\mu^{2}+\displaystyle\int_{E}\rho^{2}(e)n(de).$$
$$
\eqno{\Box}
$$
{\bf Remark 2.3.} According to the relation between It\^{o} integral
and Stratonovich integral, if we transform SDE (1.1) to the form of
Stratonovich integral, from (2.1) we know that the solution to (1.1)
enjoys SVP in a closed submanifold $K\subset R^{m}$ if and only if
the coefficients are tangent to $K$ along $K$ and the solution jumps
from $K$ to $K$. So we somewhat generalize the result in [4].

\vspace{5mm}

{\section*{\large \bf3. A special case: $S^{2}$ }}

If
$K=S^{2}=\{(x_{1},x_{2},x_{3})|x_{1}^{2}+x_{2}^{2}+x_{3}^{2}=1\}$,
then according to Theorem 2.2 we  have

{\bf Corollary 3.1. }Under the assumptions (A1) and (A2), SDE (1.1)
enjoys SVP in $S^{2}$  if and only if, $ \forall t\in [0,T],\forall
\bar{x}\in S^{2}$,
$$
\left\{\begin{array}{l} 2\langle
b(t,\bar{x}),\bar{x}\rangle+\sum\limits_{\a=1}\limits^{d}|\sigma_{\a}(t,\bar{x})|^{2}-2\displaystyle\int_{E}\langle\gamma(t,\bar{x},e),\bar{x}\rangle
n(de)=0,\\
\langle\sigma_{\a}(t,\bar{x}),\bar{x}\rangle=0,\forall
\a=1,2,...,d,\\
\bar{x}+\gamma(t,\bar{x},e)\in S^{2},n(de)-a.s..
\end{array}\right.
\eqno(3.1)
$$

\vskip 10mm In the case where there is no jump:
$$
X^{t,x}_{s}=x+\int^{s}_{t}b(r,X^{t,x}_{r})dr+\int^{s}_{t}\sigma(r,X^{t,x}_{r})dW_{r}.
\eqno (3.2)
$$

{\bf Corollary 3.2.} Under the assumptions (A1) and (A2) (without
jump), SDE (3.2) enjoys SVP in $S^{2}$ if and only if, $ \forall
t\in [0,T],\forall \bar{x}\in S^{2}$,
$$
2\langle
b(t,\bar{x}),\bar{x}\rangle+\sum\limits_{\a=1}\limits^{d}|\sigma_{\a}(t,\bar{x})|^{2}=0,
\langle\sigma_{\a}(t,\bar{x}),\bar{x}\rangle=0, \forall
\a=1,2,...,d. \eqno (3.3)
$$

In the following three examples, we set $d=1$.

{\bf Example 3.3.}
$$
\left (\begin{array}{l}
X_{s}^{1}\\
X_{s}^{2}\\
X_{s}^{3}
\end{array} \right )=
\left (\begin{array}{l}
\cos \b\\
\sin \b\\
0
\end{array} \right )+
\int^{s}_{t}\left (\begin{array}{l}
0\\
-\frac{1}{2}X_{r}^{2}\\
-\frac{1}{2}X_{r}^{3}
\end{array} \right )dr+
\int^{s}_{t}\left (\begin{array}{l}
0\\
-X_{r}^{3}\\
X_{r}^{2}
\end{array} \right )dW_{r}.
\eqno(3.4)$$

Obviously the coefficients of (3.4) satisfy (3.3), so the solution
to SDE (3.4) enjoys SVP in $S^{2}$. In fact the solution to (3.4) is
$$
X_{s}^{1}\equiv \cos \b, X_{s}^{2}=\sin \b \cos(W_{s}-W_{t}),
X_{s}^{3}=\sin \b \sin(W_{s}-W_{t}).
$$

The following counter-example shows that (3.3) is really a necessary
condition for viability.

{\bf Example 3.4.}
$$
\left (\begin{array}{l}
X_{s}^{1}\\
X_{s}^{2}\\
X_{s}^{3}
\end{array} \right )=
\left (\begin{array}{l}
\cos \b\\
\sin \b\\
0
\end{array} \right )+
\int^{s}_{t}\left (\begin{array}{l}
0\\
-\frac{3}{2}X_{r}^{2}\\
-\frac{3}{2}X_{r}^{3}
\end{array} \right )dr+
\int^{s}_{t}\left (\begin{array}{l}
0\\
-X_{r}^{3}\\
X_{r}^{2}
\end{array} \right )dW_{r}.
$$
The solution to this SDE is
$$
X_{s}^{1}\equiv \cos \b, X_{s}^{2}=\sin \b e^{-(s-t)}
\cos(W_{s}-W_{t}), X_{s}^{3}=\sin \b e^{-(s-t)}\sin(W_{s}-W_{t}).
$$
We can see easily that it can not enjoy SVP in $S^{2}$. Because in
this case,
$$
2\langle b(t,\bar{x}),\bar{x}\rangle+|\sigma(t,\bar{x})|^{2}\neq 0.
$$

{\bf Example 3.5. }Assume $\{T_{i}\}_{i=1}^{\infty}$ are i.i.d. r.v.
sequences and $T_{i}\sim \varepsilon (\lambda)$. Set
$$N_{t}:=\sup\{n\in N,\Sigma_{i=1}^{n}T_{i}\leq t\},\mbox{ \ so \ } N_{t}\sim \mbox{Poisson}(\lambda t).$$
Consider the following SDEs:
$$
\begin{array}{ll}
\left (\begin{array}{l}
X_{s}^{1}\\
X_{s}^{2}\\
X_{s}^{3}
\end{array} \right )=&
\left (\begin{array}{l}
\cos \b\\
\sin \b\\
0
\end{array} \right )+
\displaystyle\int^{s}_{t}\left (\begin{array}{l}
0\\
-\frac{1}{2}X_{r}^{2}-2\lambda X_{r}^{2}\\
-\frac{1}{2}X_{r}^{3}-2\lambda X_{r}^{3}
\end{array} \right )dr+
\displaystyle\int^{s}_{t}\left (\begin{array}{l}
0\\
-X_{r}^{3}\\
X_{r}^{2}
\end{array} \right )dW_{r}\\
&+ \displaystyle\int^{s}_{t}\displaystyle\int_{E}\left
(\begin{array}{l}
0\\
-2X_{r-}^{2}\\
-2X_{r-}^{3}
\end{array} \right )\tilde{N}(drde).
\end{array}
\eqno(3.5)
$$
where $$\int_{E}\tilde{N}(drde)=dN_{r}-n(E)dr=dN_{r}-\lambda dr.$$
Obviously the coefficients of (3.5) satisfy (3.1), so the solution
to SDE (3.5) enjoys SVP in $S^{2}$. In fact the solution to (3.5) is
$$
X_{s}^{1}\equiv \cos \b, X_{s}^{2}=\sin
\b\cos[(W_{s}-W_{t})+\pi(N_{s}-N_{t})], X_{s}^{3}=\sin \b
\sin[(W_{s}-W_{t})+\pi(N_{s}-N_{t})].
$$

\end{document}